\documentclass[12pt]{article}

\usepackage[margin=1in]{geometry}

\usepackage{graphicx,amsmath,latexsym,amssymb,cancel}
\usepackage{tikz}
\usepackage{tikz-cd}
\usetikzlibrary{arrows.meta}

\usepackage{authblk}

\usepackage{capt-of}

\usepackage{mathdots}
\usepackage[psdextra,plainpages=false,pdfpagelabels,colorlinks=false,urlcolor=black,pdfpagemode=UseNone,pdfstartview=FitH]{hyperref}
\usepackage{amsthm}
\usepackage[overload]{empheq}
\usepackage{cleveref}
\usepackage{fancyhdr}
\usepackage{float}

\newtheorem{theorem}{Theorem}[section]

\newtheorem{definition}[theorem]{Definition}

\newtheorem{lemma}[theorem]{Lemma}
\newtheorem{example}[theorem]{Example}

\newtheorem*{theorem*}{Theorem}

\newcommand{\Z}{{\mathbb Z}}

\newcommand{\sm}{\backslash}

 \newcommand{\set}[1]{\left\{#1\right\}}

\newtheorem{thm}{Theorem}[section]

\numberwithin{equation}{section}

\usepackage{tikz}
\usepackage{tikz-cd}
\usetikzlibrary{arrows.meta}
\usetikzlibrary{shapes.geometric, arrows,positioning}
\tikzstyle{rect} = [rectangle, rounded corners, 
minimum width=3cm, 
minimum height=1cm,
text centered, 
draw=black, 
fill=white]

\newcommand{\N}{{\mathbb{N}}}

\newcommand{\jjoin}{\stackrel{j}\asymp}
\newcommand{\onejoin}{\stackrel{1}\asymp}

\renewcommand{\sm}{\setminus}

\newcommand{\nfoldjoin}[2]{\underset{#1}{\stackrel{#2}{\asymp}}}

\title{Monoid Representation of the Determinants of the Adjacency Matrices of 1-joined Path Graphs}
\author{Daniel Pinzon, Daniel Pragel, and Joshua Roberts \\
Zack Marshall and Olivia Strain}
\affil{Department of Mathematics and Statistics, Georgia Gwinnett College}
\date{}

\begin{document}
\maketitle


\abstract{
    Following \cite{PPR}, this paper explores all labeled path graphs joined by 1 edge where the labeling contains the information on how to join that graph. The algebraic structure that this operation has on the determinant of the adjacency matrix of the graph gives a finite monoid of matrices under matrix multiplication. We completely characterize the monoid using the representation theory discussed in \cite{KST}. Finally, we demonstrate how these results can be used to find general formulas for the determinants of chains of 1-joined labeled paths.}

\section{Introduction}
A labeled \textit{graph} is a set of vertices, $V(G)=\{1,...,n\}$, connected by edges $E(G)$. We refer to the number of vertices in a path as its \textit{order}. A path of order $n$ is denoted as $P_n$. The \textit{standard labeling} of $P_n$ means one endpoint of the graph is labeled as 1 (the \textit{entry} vertex), then each adjacent vertex label increases by 1, resulting in the other endpoint being labeled as $n$ (the \textit{exit} vertex) for a path graph of order $n$.

To relabel a path graph $P_n$ of order $n$, we use the permutation group $S_n$.  For a path graph $P_n$, We denote $P_n|\sigma$ to be the relabeled graph where $E(P_n|\sigma) = \{\{\sigma(i),\sigma(i+1)\}\mid 1\leq i \leq n\}$  For example, $P_5|_{(1 2)(3 4 5)}$ gives the  following relabeling of $P_5$.

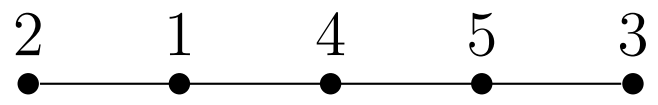
\begin{figure}[H]
    \centering
    \begin{tikzpicture}[scale=0.5]
        \tikzstyle{every node}=[circle,fill=black,inner sep=0.5mm,yshift=0,xshift=0]
        \node  [label={above:2}] (1) at  (-5,0)  {};
        \node [label={above:1}] (2) at (-3,0) {};
        \node [label={above:4}](3) at (-1,0) {};
        \node [label={above:5}](4) at (1,0) {};
        \node [label={above:3}](5) at (3,0) {};
        \draw (1)--(5);
    \end{tikzpicture}
    \caption{$G=P_5|_{(1 2)(3 4 5)}$}
    \label{fig:$P_5|_{(1 2)(3 4 5)}$}
\end{figure}

The standard labeling of $P_n$ is denoted  $P_n$ or $P_n|_{(1)}$. We define $\mathcal{P}_n=\{P_n|_{\sigma}\mid  \sigma\in S_n\}$ to be the set of all paths with any labelings. 
\bigbreak

The \textit{adjacency matrix} $A(G)$ of a graph $G$ is a square matrix in which the $(i,j)$ entry is 1 if there is an edge from the $i$th vertex of $G$ to the $j$th vertex of $G$ and 0 otherwise. We are interested in the determinants of adjacency matrices, which we will denote $|G|=\det(A(G))$. Note that relabeling a graph does not change its determinant. However, we are interested in these relabelings because they contain information on where to join graphs based on the j-join definition we will present, and the way that two graphs are joined will impact the determinant.

Some relevant proven facts about determinants of graphs are as follows:

\begin{lemma} \cite{AA}
For $n \in \N$,
$$|P_n|=\begin{cases} 
      1 & n\equiv 0\bmod 4 \\
      -1 & n\equiv 2\bmod 4 \\
      0 & n\equiv 1,3\bmod 4 
   \end{cases}
$$ 
\label{lemma:PnDet}
\end{lemma}
We note that the determinant of an isolated point $|\boldsymbol{\cdot}|=0$

\begin{lemma} \cite{AA}
If $G$ is a graph with an odd number of vertices, then $|G|$ is even. 
\label{lemma:odd_is_even_det}
\end{lemma}

\begin{lemma}
    The determinant of a graph $G$ that is a union of disjoint subgraphs $G_1,G_2,...,G_k$ is the product of the determinant of each disjoint subgraph. That is, $|G|=|G_1|\cdot|G_2|\cdot...\cdot|G_k|$ where $G$ is a disjoint union of $G_1,G_2,...,G_k$.
\label{lemma:disjoint_product_determinants}
\end{lemma} 

To chain two or more paths together, we will be using the $1$-join operation, defined in \cite{PPR}, where the labeling of each graph describes which vertices will be joined by edges. To join graph $G$ to $H$, denoted by {$G\onejoin H$}, we will form an edge from the vertex in $G$ labeled with the largest number to the vertex in $H$ labeled with 1. Then all the vertices of $G$ keep the same labels, and all of the vertices of $H$ increase their label by $|V(G)|$, as illustrated in the following example.

\begin{example} \label{1-Join Relabeling Example}
    Consider $G=P_4|_{(3 4)}$ and $, H=P_5|_{(1432)}$.
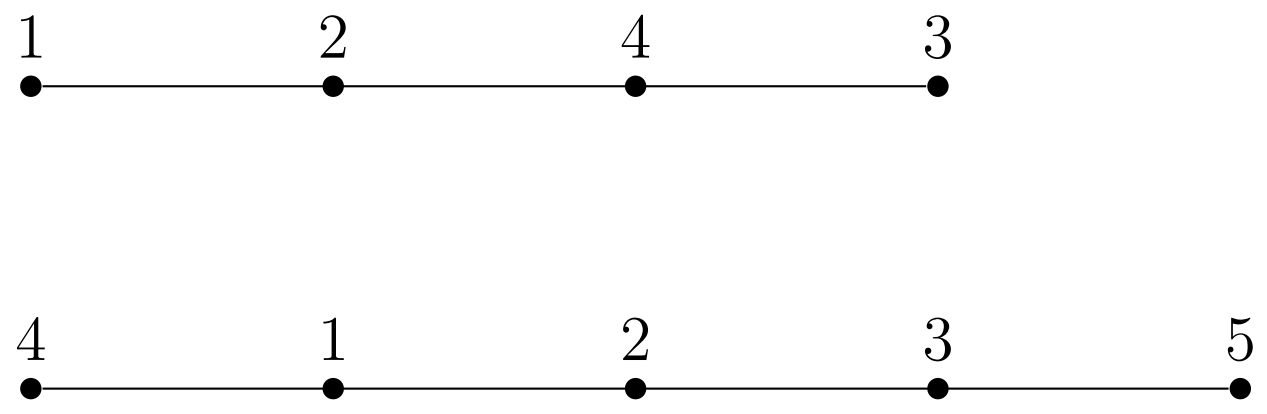
\begin{figure}[H]
\centering
\begin{tikzpicture}[scale=1]
    \tikzstyle{every node}=[circle,fill=black,inner sep=0.5mm,yshift=0,xshift=0]
    \node [label={above:1}] (1) at (-5,0)  {};
    \node [label={above:2}] (2) at (-3,0) {};
    \node [label={above:4}](4) at (-1,0) {};
    \node [label={above:3}](3) at (1,0) {};
    \draw (1)--(3);
    \node  [label={above:4}] (8) at  (-5,-2)  {};
    \node [label={above:1}] (5) at (-3,-2) {};
    \node [label={above:2}](6) at (-1,-2) {};
    \node [label={above:3}](7) at (1,-2) {};
    \node [label={above:5}](9) at (3,-2) {};
    \draw (8)--(9);
\end{tikzpicture}
\caption{$G=P_4|_{(3 4)}$, $H=P_5|_{(1234)}$}
\label{fig:pathsGH}
\end{figure}

To 1-join $G$ and $H$, we form an edge between the exit ($4$) vertex of $G$ and the entry ($1$) vertex of $H$.
We leave the labels of the vertices of $G$ the same and add $4=|G|$ to each vertex's label in $H$. 

\begin{figure}[H]
\centering
\begin{tikzpicture}[scale=1]
    \tikzstyle{every node}=[circle,fill=black,inner sep=0.5mm,yshift=0,xshift=0]
    \node [label={above:1}] (1) at (-5,0)  {};
    \node [label={above:2}] (2) at (-3,0) {};
    \node [label={above:4}](4) at (-1,0) {};
    \node [label={above:3}](3) at (1,0) {};
    \draw (1)--(3);
    \node  [label={above:8}] (8) at  (-5,-2)  {};
    \node [label={above:5}] (5) at (-3,-2) {};
    \node [label={above:6}](6) at (-1,-2) {};
    \node [label={above:7}](7) at (1,-2) {};
    \node [label={above:9}](9) at (3,-2) {};
    \draw (4)--(5);
    \draw (8)--(9);
\end{tikzpicture}
\caption{$G\onejoin H$}
\label{fig:pathsGHb}
\end{figure}
\end{example} 

We define the graph operation that we will use in this paper.
\begin{definition}[Vertex Deletion] We denote by $G\sm \{i\}$ the subgraph of $G$ obtained by removing the vertex $i$ from $V(G)$ and all edges that are incident with $i$ from $E(G).$ This is equivalent to deleting the $i$th row and column of the adjacency matrix. 
\end{definition}

\begin{definition}
    
A \textit{sandwich operation by $a$} on the set $S$, introduced in \cite{JBH} by Hickey is an operation such that given a semigroup $(S, \cdot)$ and $a\in (S, \cdot)$, we can define $x\circ y=xay$ for $x,y\in S$. Under this operation, the set $S$ is again a semigroup denoted $(S,a)$. 
\end{definition}

Following \cite{GP}, we can write the determinant of the adjacency matrix of the j-join of two graphs $G$ and $H$ as a sum of determinants of the adjacency matrices of modifications of $G$ and $H$.  The modifications made are vertex deletions of the entry (1) and exit ($n$) vertices of the graphs. Specifically for 1-joined graphs, $|G\onejoin H|=|G||H|-|G\sm\{n\}||H\sm\{1\}|$. As in \cite{PPR}, these modifications are encoded in the homomorphism from the set of graphs under the 1-join operation to $(M(2,\Z),E_1)$ the 2 by 2 integer matrices under the sandwich operation by $E_1=\begin{bmatrix}
    1&0\\0&-1
\end{bmatrix}$ defined below.

\begin{definition} Following \cite{PPR}, we define $\phi: (G_{j=1}, \onejoin) \to (M(2,\Z),E_1)$ where $G_{j=1}${are graphs where $|V(G)|\ge 2$} by $\phi(G) = \begin{bmatrix}
    |G| & |G\sm \{n\}| \\ 
    |G\sm \{1\}| & |G\sm \{1,n\}|
\end{bmatrix}$ where $G$ has $n$ vertices.
\label{def:phi}
\end{definition}

\begin{thm} \cite{PPR} The function $\phi$ defined above is a homomorphism 
\begin{equation}
    \phi(G\onejoin H)=\phi(G)E_1\phi(H)\equiv\phi(G)\circ\phi(H)
    \label{homomorphism_equation}
\end{equation}
where $E_1=\begin{bmatrix}
    1&0\\0&-1
\end{bmatrix}$.
\end{thm}

Hickey further showed that if $(S,a)$ has an identity element, then there is a semigroup isomorphism from $(S,a)$ to $(S,\cdot)$. We then note, following \cite{PPR}, that $\phi(P_4|_{(1)})=E_1$ is the identity element for $(M(2,\Z),E_1)$ where the semigroup isomorphism is defined as $\psi(M)=ME_1$. This gives the following commutative diagram.

\begin{equation}
\label{Relations_Diagram}
\begin{tikzcd}[sep=1.5cm, ampersand replacement=\&]
    \mathcal{G}_j \arrow[rd,"\pi"] \arrow[r, "\phi"] \& M(k,\Z,E_j) \arrow[d,"\psi"] \\ 
             \& M(k,\Z,\cdot) 
\end{tikzcd}
\end{equation}

\begin{example}
From \Cref{1-Join Relabeling Example}, $\phi(G\onejoin H)= \begin{bmatrix}
    0 & 1 \\
    0 & 0
\end{bmatrix}$.
This can be found directly or by computing $\phi(G\onejoin H)=
\phi(G)E_1\phi(H)=
\begin{bmatrix}
    1 & 0 \\
    0 & 0
\end{bmatrix}
\begin{bmatrix}
    1 & 0 \\
    0 & -1
\end{bmatrix}
\begin{bmatrix}
    0 & 1 \\
    0 & 0
\end{bmatrix}
=\begin{bmatrix}
    0 & 1 \\
    0 & 0
\end{bmatrix}$.
By \Cref{def:phi}, the $(1,1)$ entry of this matrix is the determinant of $G\onejoin H$.
    
\end{example}

\section{Possible Images of Path Graphs}

We now consider $\phi(P_n)$ for all labelings of $P_n$. The locations of the entry ($1$) and exit ($n$) vertices are all that impact the matrix $\phi(P_n)$ by \Cref{def:phi}. This reduces the number of labelings with possibly distinct images to $n\choose2$ 
 
\begin{example}
$P_5$ has $5!=120$ labelings but only ${5\choose2}=10$ labelings that could result in unique matrix images under $\phi$. However, from those ten, only four resulted in distinct matrices, 

$\phi(P_5|_{(1)})=\begin{bmatrix}
    0&1\\
    1&0
\end{bmatrix}
$ where $P_5|_{(1)}=$
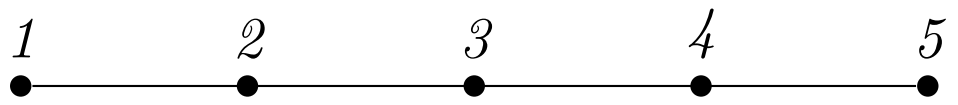
\begin{tikzpicture}[scale=.75]
\tikzstyle{every node}=[circle,fill=black,inner sep=0.5mm,yshift=0,xshift=0]
\node  [label={above:1}] (1) at  (-5,0)  {};
\node [label={above:2}] (2) at (-3,0) {};
\node [label={above:3}](3) at (-1,0) {};
\node [label={above:4}](4) at (1,0) {};
\node [label={above:5}](5) at (3,0) {};
\draw (1)--(5);
\end{tikzpicture}

$\phi(P_5|_{(12)})=\begin{bmatrix}
    0&1\\
    0&0
\end{bmatrix}$ where $P_5|_{(12)}=$
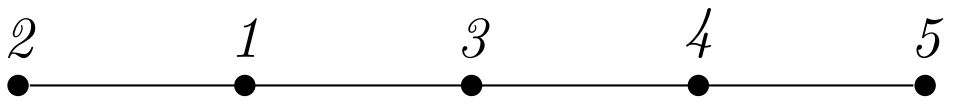
\begin{tikzpicture}[scale=.75]
\tikzstyle{every node}=[circle,fill=black,inner sep=0.5mm,yshift=0,xshift=0]
\node  [label={above:2}] (1) at  (-5,0)  {};
\node [label={above:1}] (2) at (-3,0) {};
\node [label={above:3}](3) at (-1,0) {};
\node [label={above:4}](4) at (1,0) {};
\node [label={above:5}](5) at (3,0) {};
\draw (1)--(5);
\end{tikzpicture}

$\phi(P_5|_{(125)})=\begin{bmatrix}
    0&0\\
    1&0
\end{bmatrix}$ where $P_5|_{(125)}=$
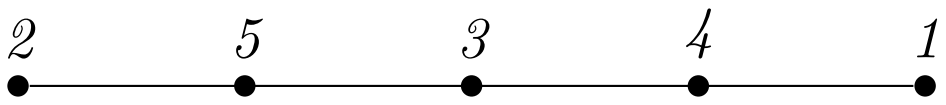
\begin{tikzpicture}[scale=.75]
\tikzstyle{every node}=[circle,fill=black,inner sep=0.5mm,yshift=0,xshift=0]
\node  [label={above:2}] (1) at  (-5,0)  {};
\node [label={above:5}] (2) at (-3,0) {};
\node [label={above:3}](3) at (-1,0) {};
\node [label={above:4}](4) at (1,0) {};
\node [label={above:1}](5) at (3,0) {};
\draw (1)--(5);
\end{tikzpicture}

$\phi(P_5|_{(12)(54)})=\begin{bmatrix}
    0&0\\
    0&0
\end{bmatrix}$
 where $P_5|_{(12)(54)}=$
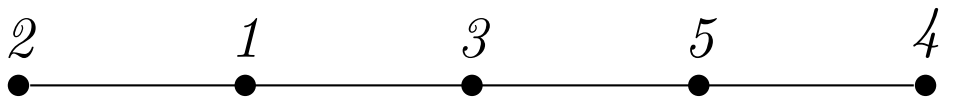
\begin{tikzpicture}[scale=.75]
\tikzstyle{every node}=[circle,fill=black,inner sep=0.5mm,yshift=0,xshift=0]
\node  [label={above:2}] (1) at  (-5,0)  {};
\node [label={above:1}] (2) at (-3,0) {};
\node [label={above:3}](3) at (-1,0) {};
\node [label={above:5}](4) at (1,0) {};
\node [label={above:4}](5) at (3,0) {};
\draw (1)--(5);
\end{tikzpicture}

\end{example}

The images of $P_n$ under $\phi$ for all possible labelings when $n=3,4,5, \text{and } 6$ give 11 distinct matrices. Before considering path orders of seven or more, we prove the following Lemma where  $\phi(\mathcal{P}_n)$ is the image under $\phi$ of all paths on $n$ vertices with all labelings.

\begin{lemma}
    The number of elements of $\phi(\mathcal{P}_n)$ is at most fifteen.
\label{theorem:15_phi_matrices}
\end{lemma}
\begin{proof}
    Recall from \Cref{def:phi} that 
    $$\phi(P_n) = \begin{bmatrix}
    |P_n| & |P_n\sm \{n\}| \\ 
    |P_n\sm \{1\}| & |P_n\sm \{1,n\}|
\end{bmatrix}$$ 
    If $n$ is odd, 
    by \Cref{lemma:PnDet}, $|\mathcal{P}_n|=0$. Removing two vertices from a path graph with an odd number of vertices results in at least one odd path, making the determinant 0. So, this reduces the number of possible matrices to nine:
    $\phi(\mathcal{P}_n)= \begin{bmatrix}
        0 & 0,\pm1 \\
        0,\pm1 & 0
    \end{bmatrix}.$ 
    
    If $n$ is even, a similar argument gives $$\phi(P_n)= \begin{bmatrix}
        \pm1 & 0 \\
        0 & 0,\pm1
    \end{bmatrix}$$
    resulting in six possible matrices.
\end{proof}

Before we prove that the 4 matrices still missing from the 15 possibilities are not the image of any labeled path graph, prove the following 2 lemmas.

\begin{lemma}
Let $G$ be a disjoint union of paths. If the order of $G$ is odd, then $|G| = 0.$
\end{lemma}

\begin{lemma}
Let $G$ and $H$ each be disjoint unions of paths of orders $m$ and $n,$ respectively. Suppose that $|G|,|H| \neq 0.$ Then, $|G| = |H|$ if and only if $m \equiv n \pmod{4}.$
\end{lemma}

\begin{proof}
By above, if $|G|,|H| \neq 0,$ then $m$ and $n$ are even, and every component of each has even order. Let $j$ and $k$ be the number of components of $G$ and $H,$ respectively, of order congruent to $2$ modulo $4.$ This implies that $|G| = (-1)^j$ and $|H| = (-1)^k.$ Additionally, we have $m \equiv 2j \pmod{4}$ and $n \equiv 2k \pmod 4.$ Thus, we see that $|G| = |H|$ if and only if $j \equiv k \pmod{2},$ which is the case if and only if $m \equiv n \pmod{4}.$
\end{proof}

\begin{theorem}\label{allPn}
     $\phi(\mathcal{P}_n)$ consists of 11 distinct matrices, specifically $\phi(P_5|_{(12)(54)}), \phi(P_5|_{(12)}),$\\
     $\phi(P_3|_{(132)}),  \phi(P_5|_{(125)}), \phi(P_3|_{(23)}) \phi(P_5|_{(1)}),  \phi(P_3|_{(1)}), \phi(P_4|_{(12)(34)}),  \phi(P_6|_{(12)(56)}),$\\ $
     \phi(P_4|_{(1)}),$ and $\phi(P_6|_{(1)})$.
\end{theorem}
\begin{proof}
The remaining 4 possible matrices are $\begin{bmatrix}
-1 & 0 \\
0 & -1 
\end{bmatrix} $,
$\begin{bmatrix}
1 & 0 \\
0 & 1 
\end{bmatrix} $, $\begin{bmatrix}
0 & -1 \\
1 & 0 
\end{bmatrix} $, and  
$\begin{bmatrix}
0 & 1 \\
-1 & 0 
\end{bmatrix}$. From the lemmas above we see looking at the non-zero entries of each matrix that none of these can be the image of a path graph. 
\end{proof}

\section{Degenerate Cases}

The $1$-join operation, defined in \cite{PPR}, requires that a graph have at least two vertices. Here we explore the image of graphs with 2 vertices or less.

By \Cref{homomorphism_equation}, for any graphs G and H, $\phi(G\onejoin P_2 \onejoin H) = \phi(G) E_1 \phi(P_2) E_1 \phi(H)$. From \cite{FH}, Harary showed that the determinant of a graph is the sum of the determinants of its possible disjoint directed cycle decompositions. Using this reasoning, $\phi(G\onejoin P_2 \onejoin H)$ has two possible disjoint directed cycle decompositions. The first would be using $P_2$ as a directed cycle and breaking down $G$ and $H$ into their disjoint directed cycle decompositions, so the product would be $|G|\cdot|H|\cdot|P_2|=-|G||H|$. The second would be forming a cycle with one vertex of $P_2$ and the exit vertex $m$ of $G$ and another cycle with the other vertex of $P_2$ with vertex 1 of $H$. This results in a product of $|G \sm \{m\}|\cdot(-1)\cdot(-1)\cdot|H\sm \{1\}|=|G\sm \{m\}||H\sm \{1\}|$. 

Thus, the (1,1) entry of $\phi(G\onejoin P_2 \onejoin H)$ is the sum $-|G||H|+|G\sm \{m\}||H\sm \{1\}|$. We assume we can define $\phi(P_2)$ in a consistent matter. Let $x$ represent the unknown value of $\phi(P_2)_{22}$. Then $\phi(P_2)=\begin{bmatrix}
    -1 & 0 \\
    0 & x
\end{bmatrix}$. By equating $\phi(G\onejoin P_2 \onejoin H)_{11}$ and $(\phi(G) E_1 \phi(P_2) E_1 \phi(H))_{11}$, we obtain the equation $-|G||H|+x|G\sm \{m\}||H\sm \{1\}| = -|G||H|+|G\sm \{m\}||H\sm \{1\}|$. Thus, $x=1$ gives a consistent result. Similar calculations allow us to define $\phi(P_2)=\begin{bmatrix}
    -1 & 0 \\
    0 & 1
\end{bmatrix}$.

Now consider a single vertex, denoted $P_1$, 1-joined with a graph $G$. Using our result above, we define $\phi(P_1)=
\begin{bmatrix}
    0 & x \\
    x & y
\end{bmatrix}
=\begin{bmatrix}
    0 & 1 \\
    1 & y
\end{bmatrix}$. Then by \Cref{def:phi}, $\phi(P_1\onejoin G)_{21}=|G|$. Equating this to $(\phi(P_1)E_1\phi(G))_{21}$, we obtain $|G|=|G|-y|G\sm \{1\}|$. Thus, $y=0$ is consistent with our results, so we can define $\phi(P_1)=\begin{bmatrix}
    0 & 1 \\
    1 & 0
\end{bmatrix}$. 

Note that $\phi(P_2)$ and $\phi(P_1)$ defined above are included in the 11 matrices established in the previous section.

\section{Monoid Representation}

The closure of the 11 distinct matrices found in \Cref{allPn} right-multiplied by $E_1$, equivalently called $\pi(\mathcal{P}_n)$ by \Cref{Relations_Diagram}, is shown in \Cref{tab:mult_table_post-E} below where we use the labeled graph to represent its matrix image under $\pi$. 

\begin{table}[H]\centering
    {\setlength{\tabcolsep}{1pt} 
    \resizebox{0.85\textwidth}{!}{
    \renewcommand{\arraystretch}{1.} 
    \begin{tabular}{c|lllccccccccll}
 & $P_5|_{(12)(54)}$& $P_4$& & & & & & & & & & &\\
  *& (0)&(id)&$P_3$&$P_3|_{(132)}$& $P_3|_{(23)}$& $P_4|_{(12)(34)}$& $P_5$& $P_5|_{(12)}$& $P_5|_{(125)}$& $P_6$& $P_6|_{(12)(56)}$& $P_5\onejoin P_5|_{(12)}$& $P_5|_{(125)}\onejoin  P_3$\\
  \hline
 $P_5|_{(12)(54)}$ (0)& $P_5|_{(12)(54)}$&$P_5|_{(12)(54)}$&$P_5|_{(12)(54)}$& $P_5|_{(12)(54)}$& $P_5|_{(12)(54)}$& $P_5|_{(12)(54)}$& $P_5|_{(12)(54)}$& $P_5|_{(12)(54)}$& $P_5|_{(12)(54)}$& $P_5|_{(12)(54)}$& $P_5|_{(12)(54)}$& $P_5|_{(12)(54)}$& $P_5|_{(12)(54)}$\\
 $P_4$ (id)& $P_5|_{(12)(54)}$&$P_4$& $P_3$& $P_3|_{(132)}$& $P_3|_{(23)}$& $P_4|_{(12)(34)}$& $P_5$& $P_5|_{(12)}$& $P_5|_{(125)}$& $P_6$& $P_6|_{(12)(56)}$& $P_5\onejoin P_5|_{(12)}$&$P_5|_{(125)}\onejoin P_3$\\
 $P_3$ & $P_5|_{(12)(54)}$&$P_3$
& $P_6$& $P_5\onejoin P_5|_{(12)}$& $P_6|_{(12)(56)}$& $P_3|_{(23)}$& $P_4$& $P_5|_{(125)}\onejoin P_3$& $P_4|_{(12)(34)}$& $P_5$& $P_5|_{(125)}$& $P_5|_{(12)}$&$P_3|_{(132)}$\\
          $P_3|_{(132)}$ & $P_5|_{(12)(54)}$&$P_3|_{(132)}$
&$P_6|_{(12)(56)}$&$P_5|_{(12)(54)}$&  $P_6|_{(12)(56)}$&  $P_5|_{(12)(54)}$&  $P_4|_{(12)(34)}$&  $P_5|_{(12)(54)}$&  $P_4|_{(12)(34)}$&  $P_5|_{(12)}$&  $P_5|_{(12)(54)}$& $P_5|_{(12)}$& $P_3|_{(132)}$\\
          $P_3|_{(23)}$ & $P_5|_{(12)(54)}$&$P_3|_{(23)}$
&$P_5\onejoin P_5|_{(12)}$&$P_5\onejoin P_5|_{(12)}$&  $P_5|_{(12)(54)}$&  $P_3|_{(23)}$&  $P_5|_{(125)}\onejoin P_3$&  $P_5|_{(125)}\onejoin P_3$&  $P_5|_{(12)(54)}$&  $P_5|_{(125)}$&  $P_5|_{(125)}$& $P_5|_{(12)(54)}$& $P_5|_{(12)(54)}$\\
          $P_4|_{(12)(34)}$ & $P_5|_{(12)(54)}$&$P_4|_{(12)(34)}$
&$P_3|_{(132)}$&$P_3|_{(132)}$&  $P_5|_{(12)(54)}$&  $P_4|_{(12)(34)}$&  $P_5|_{(12)}$&  $P_5|_{(12)}$&  $P_5|_{(12)(54)}$&  $P_6|_{(12)(56)}$&  $P_6|_{(12)(56)}$& $P_5|_{(12)(54)}$& $P_5|_{(12)(54)}$\\
          $P_5$ & $P_5|_{(12)(54)}$&$P_5$
&$P_4$&$P_5|_{(125)}\onejoin P_3$&  $P_4|_{(12)(34)}$&  $P_5|_{(125)}$&  $P_6$&  $P_5\onejoin P_5|_{(12)}$&  $P_6|_{(12)(56)}$&  $P_3$&  $P_3|_{(23)}$& $P_3|_{(132)}$& $P_5|_{(12)}$\\
          $P_5|_{(12)}$ & $P_5|_{(12)(54)}$&$P_5|_{(12)}$
&$P_4|_{(12)(34)}$&$P_5|_{(12)(54)}$&  $P_4|_{(12)(34)}$&  $P_5|_{(12)(54)}$&  $P_6|_{(12)(56)}$&  $P_5|_{(12)(54)}$&  $P_6|_{(12)(56)}$&  $P_3|_{(132)}$&  $P_5|_{(12)(54)}$& $P_3|_{(132)}$& $P_5|_{(12)}$\\
          $P_5|_{(125)}$ & $P_5|_{(12)(54)}$&$P_5|_{(125)}$
&$P_5|_{(125)}\onejoin P_3$&$P_5|_{(125)}\onejoin P_3$&  $P_5|_{(12)(54)}$&  $P_5|_{(125)}$&  $P_5\onejoin P_5|_{(12)}$&  $P_5\onejoin P_5|_{(12)}$&  $P_5|_{(12)(54)}$&  $P_3|_{(23)}$&  $P_3|_{(23)}$& $P_5|_{(12)(54)}$& $P_5|_{(12)(54)}$\\
          $P_6$ & $P_5|_{(12)(54)}$&$P_6$
&$P_5$&$P_5|_{(12)}$&  $P_5|_{(125)}$&  $P_6|_{(12)(56)}$&  $P_3$&  $P_3|_{(132)}$&  $P_3|_{(23)}$&  $P_4$&  $P_4|_{(12)(34)}$& $P_5|_{(125)}\onejoin P_3$& $P_5\onejoin P_5|_{(12)}$\\
          $P_6|_{(12)(56)}$ & $P_5|_{(12)(54)}$&$P_6|_{(12)(56)}$
&$P_5|_{(12)}$&$P_5|_{(12)}$&  $P_5|_{(12)(54)}$&  $P_6|_{(12)(56)}$&  $P_3|_{(132)}$&  $P_3|_{(132)}$&  $P_5|_{(12)(54)}$&  $P_4|_{(12)(34)}$&  $P_4|_{(12)(34)}$& $P_5|_{(12)(54)}$& $P_5|_{(12)(54)}$\\
  $P_5\onejoin P_5|_{(12)}$& $P_5|_{(12)(54)}$&$P_5\onejoin P_5|_{(12)}$&$P_5|_{(125)}$&$P_5|_{(12)(54)}$& $P_5|_{(125)}$& $P_5|_{(12)(54)}$& $P_3|_{(23)}$& $P_5|_{(12)(54)}$& $P_3|_{(23)}$& $P_5|_{(125)}\onejoin P_3$& $P_5|_{(12)(54)}$& $P_5|_{(125)}\onejoin  P_3$& $P_5\onejoin  P_5|_{(12)}$\\
  $P_5|_{(125)}\onejoin  P_3$& $P_5|_{(12)(54)}$&$P_5|_{(125)}\onejoin  P_3$&$P_3|_{(23)}$&$P_5|_{(12)(54)}$& $P_3|_{(23)}$& $P_5|_{(12)(54)}$& $P_5|_{(125)}$& $P_5|_{(12)(54)}$& $P_5|_{(125)}$& $P_5\onejoin  P_5|_{(12)}$& $P_5|_{(12)(54)}$& $P_5\onejoin  P_5|_{(12)}$& $P_5|_{(125)}\onejoin  P_3$\end{tabular}}}
    \caption{Multiplication table of the image under $\pi$ of all paths with all labelings.}
    \label{tab:mult_table_post-E}
\end{table}

This table contains two additional matrices which are not the images of path graphs, but are in the closure, namely  $P_5|_{(1)} \onejoin P_5|_{(12)}$ and  $P_5|_{(125)} \onejoin P_3|_{(1)}$. Of course, there are an infinite number of graphs that could represent any of these matrices.
\begin{figure}[H]
\begin{center}
\begin{tikzpicture}[scale=0.75]
    \tikzstyle{every node}=[circle,fill=black,inner sep=0.5mm,yshift=0,xshift=0]
    \node [label={above:1}] (1) at (-5,0)  {};
    \node [label={above:2}] (2) at (-3,0) {};
    \node [label={above:3}](3) at (-1,0) {};
    \node [label={above:4}](4) at (1,0) {};
    \node [label={above:5}](5) at (3,0) {};
    \draw (1)--(5);
    \node  [label={above:7}] (7) at  (-5,-2)  {};
    \node [label={above:6}] (6) at (-3,-2) {};
    \node [label={above:8}](8) at (-1,-2) {};
    \node [label={above:9}](9) at (1,-2) {};
    \node [label={above:10}](10) at (3,-2) {};
    \draw (5)--(6);
    \draw (7)--(10);
\end{tikzpicture}
\caption{$P_5|_{(1)} \onejoin P_5|_{(12)}$}
\label{FE graph}
\end{center}
\end{figure}
We see that we have a closed finite monoid that completely describes how the determinant changes when 1-joining  labeled path graphs. In the next two chapters we completely describe the representation of the monoid and give examples on how to easily calculate arbitrarily long chains of labeled path graphs join by 1 edge. 

\begin{theorem}
     The set of all paths with all labelings under the 1-join operation is homomorphic to a finite monoid of integer matrices under standard matrix multiplication.
\end{theorem}

\subsection{Representation Theory of a Monoid}

To better understand the structure of the monoid above, we follow Green's work in \cite{KST} on how to represent a monoid by partitioning the elements into equivalence classes. We then find which classes are isomorphic to familiar groups. Below we explain this process.
 
An element $a$ is \textit{left-less} than an element $b$ if there exists an element $c$ such that $ca=b$, denoted $a\leq_Lb$. Similarly, $a$ is \textit{right-less} than $b$ if there exists a $c$ such that $ac=b$, denoted $a\leq_Rb$. Finally, $a$ is said to be \textit{left-right-less} than $b$ if there exist $c,c'$ such that $cac'=b$, denoted $a\leq_{LR}b$.

It follows then that $a$ is \textit{left-equivalent} to $b$ if and only if $a\leq_Lb$ and $b\leq_La$. Similarly, $a$ is \textit{right-equivalent} to $b$ if and only if $a\leq_Rb$ and $b\leq_Ra$. Finally, $a$ is \textit{left-right-equivalent} to $b$ if and only if $a\leq_{LR}b$ and $b\leq_{LR}a$. The table below summarizes these definitions.

\bigbreak
\begin{minipage}{.35\textwidth}
    \begin{center}
    \underline{\textbf{Orderings}}\vspace{.15cm}\\
    $a\leq_Lb$ iff $\exists c$, $ca=b$. \\
    $a\leq_Rb$ iff $\exists c$, $ac=b$.  \\
    $a\leq_{LR}b$ iff $\exists c,c'$, $cac'=b$.
    \end{center}
\end{minipage}
\vline
\begin{minipage}{.45\textwidth}
    \begin{center}
    \underline{\textbf{Equivalence Relations}}\vspace{.15cm}\\
    $a\sim_Lb$ iff $a\leq_Lb$ and $b\leq_La$\\ 
    $a\sim_Rb$ iff $a\leq_Rb$ and $b\leq_Ra$\\ 
    $a\sim_{LR}b$ iff $a\leq_{LR}b$ and $b\leq_{LR}a$ 
    \end{center}
\end{minipage}
\bigbreak

All elements that are left-equivalent to each other are contained in a \textit{left-cell}. Similarly, all elements that are right-equivalent are contained in a \textit{right-cell}. The intersection of a left-cell and a right-cell is an equivalence class called an \textit{H-cell}. That is to say, each H-cell consists of elements of the monoid that are both left-equivalent and right-equivalent. A \textit{J-cell} contains all elements that are left-right-equivalent. A monoid that contains at least one non-invertible element $x$ can not be contained in a singular J-cell, as not all elements are left-right-equivalent to the identity element in the set.

\begin{figure}[H]
\centering
\begin{tikzpicture}[node distance=2cm]
    \node (zero) [rect] {$\set{P_5|_{(12)(54)}(0)}$};
    \node (a) [rect,below of=zero] {\oddsidemargin=0pt 
    \arraycolsep=8pt\def\arraystretch{1.5}
    $\begin{array} {c |c }
    \set{P_3|_{(123)}, P_5|_{(12)}} & \set{P_4|_{(12)(34)}, P_6|_{(12)(56)}}\\
    \hline
    \set{P_5|_{(1)} \onejoin P_5|_{(12)}, P_5|_{(125)} \onejoin P_3|_{(1)}} & \set{P_3|_{(23)}, P_5|_{(125)}}
    \end{array}$ };
        \node (one) [rect, below of=a] {$\set{P_4(1), P_3, P_5, P_6}$};
        \draw (zero)--(a);
        \draw (a)--(one);      
    \end{tikzpicture}
\caption{Diagram of the representation of the monoid}
\label{fig: monoid representation cells}
\end{figure}
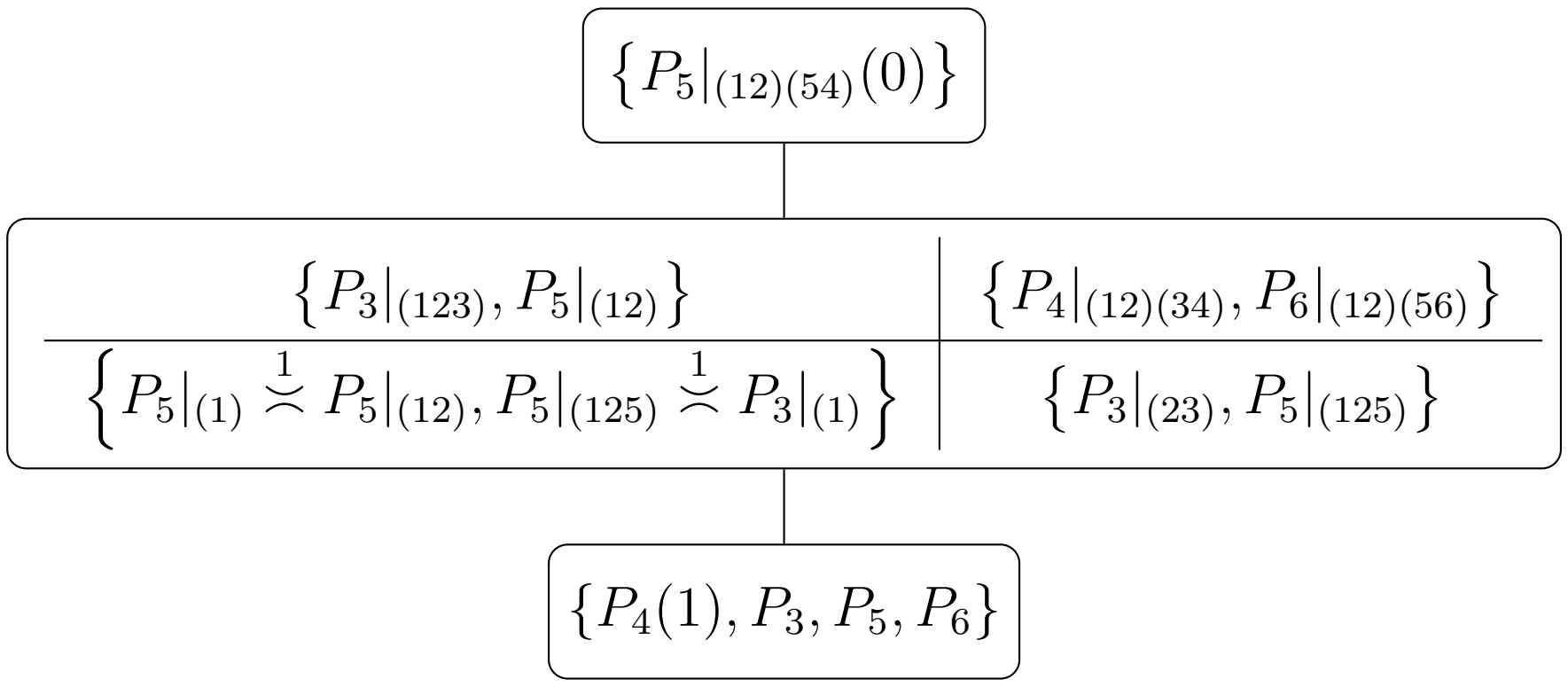

The equivalence classes of $\pi(\mathcal{P}_n)$ are shown in \Cref{fig: monoid representation cells} below where we again use the convention of denoting the matrices by a representative graph. We also denote $P_n$ to be the path graph with standard labeling.

By \cite{KST}, if an H-cell contains an idempotent, then that H-cell is isomorphic to a group, and we write $H(e)$. Since $\pi(P_4)\cdot \pi(P_4)=\pi(P_4)$, the bottom H-cell is isomorphic to $\Z_4$. Additionally, both $\pi(P_4|_{(12)(34)})\cdot \pi(P_4|_{(12)(34)})=\pi(P_4|_{(12)(34)})$ and $\pi(P_5|_{(125)}\onejoin  P_3)\cdot \pi(P_5|_{(125)}\onejoin  P_3)=\pi(P_5|_{(125)}\onejoin  P_3)$, thus $\set{P_4|_{(12)(34)}, P_6|_{(12)(56)}}$ and $\set{P_5|_{(1)} \onejoin P_5|_{(12)}, P_5|_{(125)} \onejoin P_3|_{(1)}}$ are isomorphic to $\Z_2$. We also see that $\pi(P_5|_{(12)(54)})\cdot \pi(P_5|_{(12)(54)})=\pi(P_5|_{(12)(54)})$, thus $\set{P_5|_{(12)(54)}}$ is trivially isomorphic to $\Z_1$.

To summarize, this monoid representation fully characterizes how the determinant of a collection of 1-joined path graphs is affected by 1-joining it to another path.

\section{Applications}
In this section, we will examine how these results can be used to find the determinants of large chains of labeled path graphs. First, note the (1,1) entry of $\pi(G)$ is equal to the (1,1) entry of $\phi(G)$.

 
A $k$-chain of a graph $G$ is the $j$-join of $k$ graphs $G$
$$ 
\nfoldjoin{k}{j} G= \underbrace{G \jjoin G \jjoin \cdots \jjoin G}_{k \text{ copies of }G}.
$$

Since $\pi$ is a homomorphism,  $\pi\left(\nfoldjoin{k}{j} G\right) =\underbrace{\pi(G) \cdots \pi(G)}_{k \text{ copies}}$ under matrix multiplication. Thus, to find the determinant of any chain of 1-joined labeled path graphs, we find the relevant products in our multiplication table and examine the (1,1) entry of the resultant matrix.


\begin{example}
    To illustrate, choose $j=1$, $k=2025$, and $G=P_3|_{(23)}\onejoin P_5|_{(125)}$. This would be the graph of a 1-joined chain of alternating $P_3|_{(23)}$ and $P_5|_{(125)}$ graphs, with a total number of 5050 path graphs being 1-joined together. Then by \Cref{tab:mult_table_post-E}, $\pi\left(\nfoldjoin{2025}{1} G\right)=(\pi(P_3|_{(23)})\cdot \pi(P_5|_{(125)}))^{2025}=(\pi(P_5|_{(12)(54)}))^{2025}= \begin{bmatrix}
    0&0\\
    0&0
\end{bmatrix}$. Thus, the determinant of this chain is 0.
\end{example}

\begin{example}
   As a final example, we choose $j=1$, $k=n$ with $n\in\N$, and $G=P_4|_{(12)(34)}\onejoin P_6|_{(12)(56)}$. This is a graph of a 1-joined chain of alternating $P_4|_{(12)(34)}$ and $P_6|_{(12)(56)}$ graphs, with a total number of $n$ path graphs 1-joined together as shown in the figure below.

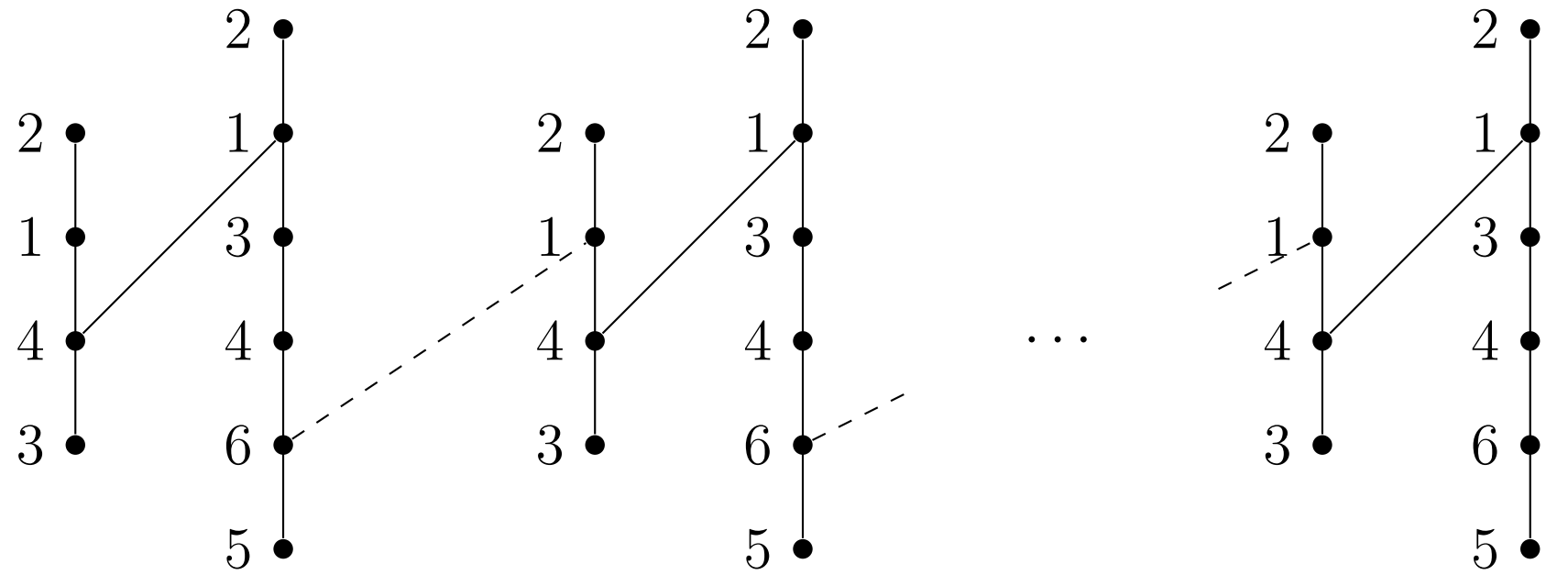
\begin{figure}[H]
\centering
    \begin{tikzpicture}[scale=0.75]
        \tikzstyle{every node}=[circle,fill=black,inner sep=0.5mm,yshift=0,xshift=0]
        \foreach \x in {-8,-3,4}{
            \node  [label={[xshift=0cm ]left:2}] (2+\x) at  (\x,4)  {};
            \node [label={[xshift=0cm ]left:1}](1+\x) at (\x,3) {};
            \node [label={[xshift=0cm ]left:4}](4+\x) at (\x,2) {};
            \node [label={[xshift=0cm ]left:3}] (3+\x) at (\x,1) {};
            \draw (2+\x)--(3+\x);
        }
        \foreach \x in {-6,-1,6}{
            \node  [label={[xshift=0cm ]left:2}] (2b+\x) at  (\x,5)  {};
            \node [label={[xshift=0cm ]left:1}](1b+\x) at (\x,4) {};
            \node [label={[xshift=0cm ]left:3}](3b+\x) at (\x,3) {};
            \node [label={[xshift=0cm ]left:4}] (4b+\x) at (\x,2) {};
            \node [label={[xshift=0cm ]left:6}](6b+\x) at (\x,1) {};
            \node [label={[xshift=0cm ]left:5}] (5b+\x) at (\x,0) {};
            \draw (2b+\x)--(5b+\x);
        }
        \draw (4+-8)--(1b+-6);
        \draw (4+-3)--(1b+-1);
        \draw (4+4)--(1b+6);
        \draw[dashed](6b+-6)--(1+-3); 
        \draw[dashed](6b+-1)--(0,1.5); 
        \draw[dashed](3,2.5)--(1+4); 
        \tikzstyle{every node}=[circle,fill=none,inner sep=0.5mm,yshift=0,xshift=0]\node[draw=none] (ell) at (1.5,2) {$\cdots$};
        \end{tikzpicture}
\caption{Chain of $P_4|_{(12)(34)}\nfoldjoin{n}{1} P_6|_{(12)(56)}$}
    
\end{figure}

Thus $\pi\left(\nfoldjoin{n}{1} G\right)=(\pi(P_4|_{(12)(34)})\cdot \pi(P_6|_{(12)(56)}))^n=(\pi(P_6|_{(12)(56)}))^n$. We see from \Cref{tab:mult_table_post-E} that $\pi(P_6|_{(12)(56)})\cdot \pi(P_6|_{(12)(56)})=\pi(P_4|_{(12)(34)})$ and $\pi(P_4|_{(12)(34)})\cdot \pi(P_6|_{(12)(56)})=\pi(P_6|_{(12)(56)})$.
\bigbreak Therefore, $\pi\left(\nfoldjoin{n}{1} G\right)=
\begin{cases}
    \pi(P_6|_{(12)(56)})=\begin{bmatrix}
    -1&0\\
    0&0
\end{bmatrix}, & n$ is odd$ \\
    \pi(P_4|_{(12)(34)})=\begin{bmatrix}
    1&0\\
    0&0
\end{bmatrix}, & n$ is even$
\end{cases}$. \bigbreak
This shows  $\left|\nfoldjoin{n}{1} G\right|=(-1)^n$. 
\end{example}

\section{Future Work}
As our work was focused strictly on 1-joining labeled path graphs, we are interested in the much more difficult task of expanding these results to 2-joins or j-joins in general. Additionally, we are interested in expanding this work to include trees and joinings of trees instead of paths. Which trees can be formed by one-joining paths? Can we classify those trees?

\bibliographystyle{unsrt}

\begin{thebibliography}{00}
\bibitem[1]{BN} A. Bar-Noy; J. Naor, \emph{Sorting, Minimal Feedback Sets and Hamilton Paths in Tournaments}, SIAM Journal of Discrete Mathematics, \textbf{3}(1) (1990), 7 - 20.
\bibitem[2]{LG} L. Goldberg, \emph{Matrix Theory with Applications}, McGraw-Hill International Editions, Mathematics and Statistics Series, 1991.
\bibitem[3]{GP} B. Gyurov, D. Pinzon, \emph{The Determinant of Graphs Joined by $j$ edges}, Thai Journal of Mathematics, \textbf{13} (2016) 2: 357-367.
\bibitem[4]{FH} F. Harary, \emph{The Determinant of the adjacency matrix of a graph}, SIAM Rev., \textbf{4} (1961), 202 - 210.
\bibitem[5]{JBH} J.B. Hickey, \emph{Semigroups Under A Sandwich Operation}, Proceedings of The Edinburgh Mathematical Society, \textbf{26(3)} (1983) 371-382, url: \url{https://doi.org/10.1017/S0013091500004442}.
\bibitem[6]{HMR} H.M. Rara,\emph{Reduction procedures for calculating the determinant of the adjacency matrix of some graphs and the singularity of square planar grids}, Discrete Mathematics, \textbf{151} (1996), 213-219, url: \url{https://core.ac.uk/download/pdf/81194617.pdf}.
\bibitem[7]{SAG}S. Sookyang, S. Arworn, B. Gyurov, \emph{Determinant of Graphs joined by two edges}, Thai Journal of Mathematics, \textbf{10} (2012) 1: 101-111.
\bibitem[8]{PPR}D.Pinzon, D. Pragel, J. Roberts (2024). Algebraic Structures on Graphs Joined by Edges, to appear in Discussiones Mathematicae - General Algebra and Applications. arXiv preprint arXiv:2409.02355.
\bibitem[9]{AA}A. Alireza, \emph{Determinants of Adjacency Matrices of Graphs}, Transactions on Combinatorics, \textbf{1} (2012) 4: 9-16 url: \url {https://toc.ui.ac.ir/article_2041.html}
\bibitem[10]{KST}M. Khovanov, M. Sitaraman, D. Tubbenhauer, \emph{Monoidal categories, representation gap and cryptography}, Transactions of the American Mathematical Society, Series B, \textbf{11} (2024) 10: 329-395 url: \url{http://dx.doi.org/10.1090/btran/151}

\end{thebibliography}

\end{document}